\documentclass[twoside]{article}

\usepackage{amsmath,amssymb,latexsym}
\usepackage{mathrsfs}
\usepackage{amsthm} 
\usepackage{color}

 \newtheorem{thm}{Theorem}[section]
 \newtheorem{cor}[thm]{Corollary}
 \newtheorem{lem}[thm]{Lemma}
 \newtheorem{prop}[thm]{Proposition}    
 \newtheorem{ex}[thm]{Example}
\theoremstyle{definition} 
\newtheorem{rem}[thm]{Remark}
 \newtheorem{defn}[thm]{Definition}
 \numberwithin{equation}{section}

\newcommand{\h}{{\mathcal H}}
\newcommand{\dom}{{\mathcal D}}
\newcommand{\ran}{{\mathcal R}}

\newcommand{\mn}{\mathbb N}

\newcommand{\mg}{\mathbb{G}}
\newcommand{\mf}{\mathbb{F}}

\newcommand{\mw}{\mathbb{W}}

\newcommand{\seqgr[1]}{(#1_i)_{i=1}^\infty }

\newcommand{\newin}{\kern-0.3em\in\kern-0.2em}
\newcommand{\newsubset}{\kern-0.2em\subset\kern-0.2em}

\newcommand{\sumii}{\sum_{i=1}^{\infty}}
\newcommand{\range}{{\mathcal R}}

\newcommand{\bp}{\noindent{\bf Proof. \ }}
\newcommand{\ep}{\noindent{$\Box$}}

\newcommand{\lin}{{\rm span}}

\title{Characterization of atomic decompositions, Banach frames, Xd-frames, duals and synthesis-pseudo-duals, with application to Hilbert frame theory}

\author{Diana T. Stoeva\\ \ \\
Acoustics Research Institute\\ 
Austrian Academy of Sciences \\
Wohllebengasse 12-14, Vienna 1040, Austria \\
dstoeva@kfs.oeaw.ac.at}

\begin{document}

\maketitle

\thispagestyle{empty}

\begin{abstract}
In this paper we consider series expansions via a frame and a non-frame and the possibilities for interchange of the two sequences, both in the Hilbert and Banach space setting. First we give a characterization of frame-related concepts in Banach spaces (atomic decompositions, Banach frames, $X_d$-Riesz bases, $X_d$-frames, $X_d$-Bessel sequences, and sequences satisfying the  lower  $X_d$-frame  condition). 
We also determine necessary and sufficient conditions for operators to preserve the type of the 
concepts listed above. 
Then we discuss differences and relationships between expansions in a Banach space 
and its dual space when interchanging the involved sequences. 
Finally, we apply some of the results to answer problems in Hilbert frame theory.  
We show that interchanging a frame and a non-Bessel sequence in series expansions is not always possible
(leading to differentiation of  analysis- and synthesis-pseudo-duals of a frame) 
and determine an appropriate subspace where interchange can be done. 
 We characterize all the synthesis-pseudo-duals of a frame and determine a class of frames whose synthesis-pseudo-duals (resp. analysis-pseudo-duals) are necessarily frames. We also investigate connections between the lower frame condition and series expansions. 
Examples are given to illustrate statements in the paper and to show the optimality of some results.  
\end{abstract}


\vspace{.05in}

\section{Introduction and Basic Definitions}
The frame-concept was introduced by Duffin and Schaeffer \cite{DSframe} in 1952.
The sequence $\seqgr[g]$ is called a {\it (Hilbert)
 frame for the Hilbert space $\h$ with bounds $A,B$} if $A$ and $B$ are positive constants and $A \|h\|^2\leq \sumii | \langle  h,g_i\rangle|^2 \leq
B\|h\|^2$ for every $h\in\h$. 
It took
several decades for scientists to realize the high potential 
 of frames. Around 1990, the frame-theory began to develop in connection with Gabor analysis and wavelets \cite{D90, Dbook, DGM}.
  Nowadays, frames are very important both for theory and real life. 
  They play fundamental role in signal and image processing and 
  find applications in  
 wireless communication, speech recognition,  geophysics,  biology,  
and many other areas. 
  For more on frame theory we refer to \cite{CasArt, Cbook, Gbook, Hbook}.
  
What makes frames very useful is that they require less restrictive conditions on the sequence elements compared to orthonormal bases and still they allow reconstructions of the elements of the Hilbert space: if $\seqgr[g]$ is a frame for $\h$, then there exists a frame $\seqgr[f]$ for $\h$ so that
\begin{equation} \label{fexp2}
f=\sum_{i=1}^\infty \langle  f,g_i \rangle f_i, \ \forall f\in\h,
\end{equation}
and
\begin{equation} \label{fexp1}
g=\sum_{i=1}^\infty \langle  g,f_i \rangle g_i, \ \forall g\in\h;
\end{equation}
 such a frame $\seqgr[f]$ is called a {\it dual frame of $\seqgr[g]$}. 
  When a frame $\seqgr[g]$  for $\h$ is at  the same time a Schauder basis of $\h$ 
  (so called {\it Riesz basis\footnote{Riesz bases were introduced by Bari \cite{Bary46,Bari} in a different but equivalent way; for more on Riesz bases and equivalent definitions we refer to \cite{Cbook,Young}.} for $\h$})
  there is only one sequence $\seqgr[f]$ satisfying (\ref{fexp2}),  resp. (\ref{fexp1}), 
 and it is also a Riesz basis for $\h$. 
   When a frame   for $\h$ is not a Schauder basis for $\h$ (so called {\it overcomplete frame for $\h$}), there are many frames $\seqgr[f]$ for $\h$ satisfying 
 (\ref{fexp2}) and (\ref{fexp1}), 
 and this property makes the overcomplete frames very attractive for applications. 
 For example, it gives the possibility to search for dual frames fulfilling some additional requirements, 
see e.g. the case of wavelet frames \cite{BL,DH}.

 The dual frames might not be the only sequences giving series expansions. 
For some overcomplete frames $\seqgr[g]$ for $\h$, in addition to the dual frames (which always satisfy both (\ref{fexp2}) and (\ref{fexp1})), there exist non-frame sequences 
 $\seqgr[f]$  satisfying both
 (\ref{fexp2}) and (\ref{fexp1}) (see Example \ref{adn} with $p=2$) 
 or satisfying only one of (\ref{fexp2}) and (\ref{fexp1}) (see Example \ref{eximp});  
 in the wavelet setting, an example of a frame $\seqgr[g]$ and a non-frame $\seqgr[f]$ satisfying (\ref{fexp2})  
 can be found in \cite{LO}. 
 This motivates the investigation of sequences $\seqgr[f]$ which are not necessarily frames, but satisfy (\ref{fexp2}) and/or (\ref{fexp1}),  
 and we naturally refer to them as {\it dual sequences of $\seqgr[g]$}. 
   Such sequences can be important for numerical stability of representations of type (\ref{fexp1}), 
    because although they might not be frames, they would necessarily satisfy the lower frame condition \cite{CCLL,LO}.  
 Furthermore, they turned out to be also involved in representations of the inverse of a frame multiplier \cite{SBsampta15} which may further lead to their use in areas where multipliers are applied, e.g. in sound synthesis \cite{DMT}, psychoacoustical modeling \cite{BLED},  denoising \cite{MBKD}.

Example \ref{eximp} motivates separate investigation of sequences, satisfying (\ref{fexp2}) or (\ref{fexp1}). 
 Given a frame $\seqgr[g]$ for $\h$, a sequence $\seqgr[f]$ with elements in $\h$ 
will be called 

- a {\it synthesis-pseudo-dual} (in short, {\it s-pseudo-dual}) {\it of  $\seqgr[g]$}, if it  satisfies (\ref{fexp2});

- an {\it analysis-pseudo-dual} (in short, {\it a-pseudo-dual}) {\it of  $\seqgr[g]$}, if it  satisfies (\ref{fexp1}).

While characterizations of all the dual frames of a given frame are well known in the literature, see e.g. \cite{CHL, Cbook, HL, Li}, less is known in the direction of non-frame dual sequences.  
In the present paper (in Section \ref{s5}) 
we investigate $s$- and $a$-pseudo-duals of frames. 
Note that throughout the paper, the series representations in (\ref{fexp2}) and (\ref{fexp1}) are always 
considered in the sense of norm convergence.
Investigation of such representations in a weak sense, i.e., representations in the form 
$\langle f,g\rangle = \sum_{i=1}^\infty \langle  f,g_i \rangle \langle  f_i, g\rangle , \ \forall f,g\in\h$ (called pseudo-frame representations),
 is done in \cite{LO};  
just to be noticed that the statement of \cite[Prop. 4.10]{LO} is not correct, 
 see  Example \ref{eximp} and the comments after it. 
For investigation of series expansions in a subspace $\h_0$ of $\h$ via a frame for $\h_0$ and a sequence 
with elements not necessarily in $\h_0$, we refer to \cite{FW}.

A natural extension of the frame inequalities to Banach spaces leads to the concepts of {\it $p$-frame} \cite{AST} and {\it $X_d$-frame} \cite{CCS} (see Definition \ref{def1}). In contrast to the frame-case, the $X_d$-frame inequalities do not necessarily lead to reconstruction via
series expansions, see \cite[Ex. 2.8]{CCS}. 
With aim to have reconstructions, {\it atomic decompositions} (giving reconstruction via series expansions) were considered by Feichtinger and Gr\"{o}chenig   \cite{FG, FG2, G} and the concept of {\it Banach frame} (giving reconstruction via an operator) was introduced in \cite{G}.  
Historically, first atomic decompositions and Banach frames were introduced, and after that the concepts of  $p$-frame and  $X_d$-frame appeared.
Further, frames were extended to Fr\`echet spaces \cite{PS2,PS4}, but for the purpoces of the paper we stay in the context of Banach spaces. 

\begin{defn} \label{def1}
Let $X$ be a Banach space, $X_d$ be a $BK$-space (i.e., a Banach sequence space for which the coordinate functionals are continuous) and \mbox{$g_i\in X^*$,} $i\in\mn$. 
The sequence $\seqgr[g]$ is called a {\it Banach frame for $X$ with respect to  $X_d$} if 
\begin{itemize}
\item[{\rm (i)}] $(g_i(f))_{i=1}^\infty \in X_d$, $\forall f\in X$;
\item[{\rm (ii)}] $\exists$ positive constants $A$ and $B$ so that $ A\|f\|_X \leq \|(g_i(f))_{i=1}^\infty\|_{X_d} \leq B \|f\|_X$, $\forall f\in X$; 
\item[{\rm (iii)}] $\exists$ bounded operator  $Q: X_d \to X$ so that $Q(g_i(f))_{i=1}^\infty =f$,  $\forall f\in X$. 
\end{itemize}
Such operator $Q$ is called a {\it Banach frame operator of $\seqgr[g]$}. 

\vspace{.05in}
Let $f_i\in X$, $i\in\mn$. The pair $(\seqgr[g], \seqgr[f])$ is called an {\it atomic decomposition  of $X$ with respect to $X_d$} if (i) and (ii) hold 
and
\begin{itemize}
\item[{\rm (iii')}]  $f=\sum_{i=1}^\infty g_i(f) f_i,  \ \forall f\in X$.
\end{itemize}

The sequence  $\seqgr[g]$ is called an {\it $X_d$-frame for $X$} 
(resp. {\it $X_d$-Bessel sequence for $X$})  if (i) and (ii) (resp (i) and the upper inequality in (ii)) hold. 

It is said that $\seqgr[g]$ satisfies the {\it lower $X_d$-frame condition} if the lower inequality in (ii) holds for all those $f$ for which $(g_i(f))_{i=1}^\infty\in X_d$.

When $X_d=\ell^p$, $p\in(1,\infty)$, then {\it $p$-frame}, {\it $p$-Bessel sequence}, and {\it lower $p$-frame condition} stand for $\ell^p$-frame, $\ell^p$-Bessel sequence, and lower $\ell^p$-frame condition, respectively; in the case when $p=2$ and $X$ is a Hilbert space, these 
are the well established concepts of a {\it frame}, {\it Bessel sequence}, and {\it lower frame condition}, respectively. 
\end{defn}

In the present paper (in Section \ref{s3}) we characterize all the concepts from the above definition. 
 When $X_d=\ell^2$ and $X$ is a Hilbert space, characterizations of the corresponding concepts can be found in \cite{BSA, Cbook} and references therein. Note that if $X_d=\ell^2$ and $X$ is a Hilbert space, then $X_d$-frame, Banach frame and the first sequence in an atomic decomposition pair mean the same (namely, a Hilbert frame). 
For general Banach spaces, these three types of sequences 
do not mean the same; for a detail discussion about their relationship and differences see \cite{Spert}.

 Riesz bases were also generalized to Banach spaces, under the name of $X_d$-Riesz bases. The concept of {\it $X_d$-Riesz basis} was established by Feichtinger and Zimmermann \cite{FZim}. Another definition for an $X_d$-Riesz basis (Definition \ref{def2} below) is considered in \cite{Srbasis}, motivated by the definitions of a {\it $p$-Riesz basis} in \cite{AST,CS}. When $X_d$  has the canonical vectors as a Schauder basis, the definitions in \cite{FZim} and \cite{Srbasis} are equivalent \cite[Sec. 3]{Srbasis}. 
\begin{defn} \label{def2} Let $Y$ be a Banach space, $X_d$ be a Banach sequence space, and $g_i\in Y$, $i\in\mn$.   
The sequence $\seqgr[g]$ is called an {\it $X_d$-Riesz basis for $Y$ with bounds A,B}, if it is complete in $Y$, the constants $A$ and $B$ are positive, and
\begin{eqnarray} \label{bo}
A \left\|\seqgr[c]\right\|_{X_d} \le \left\| \sum_{i=1}^\infty c_i g_i \right\|_{Y} \le
B\left\|\seqgr[c]\right\|_{X_d}, \ \forall \seqgr[c]\in X_d.
\end{eqnarray}
\end{defn}

The paper is organized as follows. 

Section \ref{notprel} concerns the notation used in the paper and some preliminaries. 
In Section \ref{s3} we give a characterization of  the sequences defined in \ref{def1} and \ref{def2} based on an operator acting on the canonical basis of the corresponding sequence space. 
We also consider operators which preserve the sequence type; we determine necessary and sufficient conditions, discussing the case of bounded operators and the case of not necessarily bounded ones.
Some of the results in Section  \ref{s3} are needed for the main part of the paper (Sections \ref{s4} and \ref{s5}), while others are included for the sake of having complete list with characterization of the popular frame-related concepts in Banach spaces,  
which is of independent interest as well. 

In Section \ref{s4} (the Banach space setting) and Section \ref{s5} (the Hilbert space setting) we consider the topic of 
interchange of sequences in series expansions, focusing on expansions involving a \lq\lq frame\rq\rq\ and a dual sequence not necessarily being a \lq\lq frame\rq\rq. 
In Section \ref{s4},  we discuss connections between expansions in $X$ and $X^*$ via an $X_d$-frame $\mg$ for $X$ and a sequence $\mf$, namely, representations of the form $f=\sum_{i=1}^\infty g_i(f) f_i, \,  f\in X,$ and representations of the form $g=\sum_{i=1}^\infty g(f_i) g_i,  \, g\in X^*$. 
Given an $X_d$-frame $\mg$ for $X$, we characterize all the sequences $\mf$ which give atomic decompositions $(\mg, \mf)$ of $X$ with respect to $X_d$;  
among those $\mf$, we characterize the ones which have \lq\lq dual\rq\rq\ properties, namely, which are $X_d^*$-frames for $X^*$. 
In Section \ref{s5} we focus on the Hilbert space setting as being of utmost importance for applications, and we
solve some problems in frame theory. The section contains not only consequences of results from Section \ref{s4}, but more extensive investigation. 
Via Example \ref{eximp} we show that  the interchange of a frame and a non-Bessel sequence in series expansions is not always possible and this motivates the separate investigation of a- and s-pseudo-duals of a frame. 
We determine an appropriate subset of the space where an interchange can be done and an example shows the optimality of the determined subset (though, in particular cases one can have a bigger subset). 
We also characterize all the s-pseudo-duals of a frame. Given a sequence $\seqgr[f]$ which satisfies the lower frame condition, we investigate the existence of representations in the form (\ref{fexp2}) via a Bessel sequence $\seqgr[g]$. Finally, we determine a class of frames whose s-pseudo-duals and a-pseudo-duals are necessarily frames.

\section{Notation and Preliminaries} \label{notprel}

Throughout the paper, $X$ denotes an infinite-dimensional separable Banach space and $X^*$ denotes its dual; 
$X_d$ denotes a Banach sequence space and $X_d^*$ denotes its dual; $\h$ denotes an infinite-dimensional separable Hilbert space and $\seqgr[e]$ denotes an orthonormal basis for $\h$. Unless otherwise specified, the letter $\mg$ (resp. $\mf$) means a sequence $\seqgr[g]$ (resp. $\seqgr[f]$) with elements from $X^*$ (resp. $X$).  
For convenience of the writings we use $\mn$ as an index set, but any countably infinite index set can be used instead. 
A linear mapping is called an {\it operator}.  
The domain (resp., the range) of an operator $V$ is denoted by $\dom(V)$ (resp., $\ran(V)$). 
The space $X_d$ is called a {\it $BK$-space} if the coordinate functionals are continuous. If the canonical vectors form a Schauder basis for $X_d$, then $X_d$ is called a {\it $CB$-space} and the canonical basis is denoted by $\seqgr[\delta]$. 
In particular, the canonical basis of $\ell^2$ is denoted by  $\seqgr[\delta]$.
When $X_d$ is a $CB$-space, then $X_d^\circledast \mathrel{\mathop:}=\{
(V \delta_i)_{i=1}^\infty : V\in X_d^* \}$ with the norm 
$\left\|(V \delta_i)_{i=1}^\infty\right\|_{X_d^\circledast}\mathrel{\mathop:}=\left\|V\right\|_{X_d^*}$ is a $BK$-space isometrically
isomorphic to $X_d^*$ 
{\rm \cite[Ch.VI \S 1]{KA59}} 
and for the rest of the paper $X_d^*$ is identified with $X_d^\circledast$; 
furthermore, for every $\seqgr[c]\in X_d$ and every $\seqgr[d]\in X_d^*$, the series $\sum_{i=1}^\infty c_i d_i$ converges and $$|\sum_{i=1}^\infty c_i d_i|\leq \|\seqgr[c]\|_{X_d}\cdot \|\seqgr[d]\|_{X_d^*}. $$
If $X_d$ is both reflexive and a $CB$-space (called an {\it $RCB$-space}), 
 then $X_d^\circledast$ is a $CB$-space and its canonical basis is denoted by $(\delta_i^*)_{i=1}^\infty$. 
The set of all finite linear combinations of elements from $\{x_i\}_{i=1}^\infty$  is denoted by 
$\lin\{x_i\}$.
The abbreviation \lq wrt\rq \, stands for \lq with respect to\rq.

For given $CB$-space $X_d$ and given $\mg$, 
being an $X_d$-Bessel sequence for $X$ or satisfying the lower $X_d$-frame condition,
the {\it analysis operator $U_\mg$} 
 and the {\it synthesis operator $T_\mg$} are determined by 
 \begin{eqnarray}
 U_\mg:\dom(U_\mg)\to X_d, & & U_\mg f = (g_i(f))_{i=1}^\infty, 
 \label{opu} \\
 T_\mg: \dom(T_\mg) \to X^*, & & T_\mg\seqgr[c]=\sum_{i=1}^\infty c_i g_i, 
 \label{opt} 
 \end{eqnarray} 
 where $\dom(U_\mg)=\{ f\in X : (g_i(f))_{i=1}^\infty\in X_d\}$ and $\dom(T_\mg)=\{ \seqgr[c]\in X_d^* : \sum_{i=1}^\infty c_i g_i \mbox{ converges in } X^*\}$.

\sloppy
Note that the definition of an $X_d$-Bessel sequence $\mg$ for $X$ requires $\dom(U_\mg)=X$, while the definition of a sequence satisfying the lower $X_d$-frame condition allows the domain of its analysis operator to be a subset of $X$.

If $X_d$ is an $RCB$-space and $\mg$ is an $X_d$-Bessel sequence for $X$, then $\dom(T_\mg)=X_d^*$,  $U_\mg^*=T_\mg$, and $T_\mg^*=U_\mg$ \cite{CCS, Sduals}.
If $X_d$ is a $BK$-space and $\mg$ satisfies the lower $X_d$-frame condition, 
then  $\ran(U_\mg)$ is closed in $X_d$ and $U_\mg$ has a bounded inverse $U_\mg^{-1}: \ran(U_\mg)\to \dom(U_\mg)$ \cite{Sduals}.

\vspace{.1in}
If $(g_i)_{i=1}^\infty$ is a {\it Bessel sequence in $\h$} with bound $B$ (i.e., $\sumii | \langle  h,g_i\rangle|^2 \leq
B\|h\|^2$, $\forall h\in\h$) and $(f_i)_{i=1}^\infty$ satisfies (\ref{fexp1}), 
then $(f_i)_{i=1}^\infty$ satisfies the {\it lower frame condition} with bound $1/B$ (i.e., $\frac{1}{B} \|h\|^2\leq \sumii | \langle  h,g_i\rangle|^2 $ for all those $h\in\h$ for which $\sumii | \langle  h,g_i\rangle|^2<\infty$) \cite{CCLL}. 
The extension of this statement to Banach spaces is 
given below and 
it can be proved easily using some calculations from the proof of \cite[Prop. 3.4]{CCS}.

\begin{lem} \label{lc}
Assume that $X_d$ is a $CB$-space. 
Let $\mg$ be an $X_d$-Bessel sequence for $X$ with bound $B$ and let there exist $\mf$ so that $f=\sum_{i=1}^\infty g_i(f) f_i$ for every $f\in X$ or $g=\sum_{i=1}^\infty g(f_i) g_i$ for every $g\in X^*$. Then $\mf$ satisfies the lower $X_d^*$-frame condition  with bound $1/B$, i.e., $\frac{1}{B}\|g\|_{X^*} \leq \|(g(f_i))_{i=1}^\infty\|_{X_d^*}$ for those $g\in X^*$ for which  $(g(f_i))_{i=1}^\infty\in X_d^*$.
\end{lem}

\section{Characterization of frame-related concepts in Banach spaces} \label{s3}

In this section we characterize the frame-related concepts from Definitions \ref{def1} and \ref{def2}, and consider operators which preserve the sequence type. 
As mentioned in the Introduction, some of the results in this Section are needed for the main part of the paper (Sections \ref{s4} and \ref{s5}), while others are included for the sake of completeness and are of independent interest. 
Notice that the characterizations concerning $X_d$-Bessel sequences, $X_d$-Riesz bases, and sequences satisfying the lower $X_d$-frame condition extend naturally results for the corresponding Hilbert space concepts \cite{Cbook,BSA}. 
The situation with $X_d$-frames, Banach frames, and atomic decompositions, is that these are three different types of extension of Hilbert frames.  
Note that the roles of $\mg$ and $\mf$ in an atomic decomposition ($\mg, \mf$) are not symmetric, so we also characterize each one of them.

\begin{thm} \label{allfr} Let $X_d$ be an $RCB$-space and let $X$ be reflexive. 
\begin{itemize} 
\item[{\rm (i)}] The $X_d$-Bessel sequences for $X$ are precisely 
the sequences $(T\delta_i^*)_{i=1}^\infty$,
where $T:X_d^*\to X^*$ is a bounded operator.

\item[{\rm (ii)}] The sequences $G$ 
satisfying the lower $X_d$-frame condition 
are precisely the sequences $(T \delta_i^*)_{i=1}^\infty$,
where the operator  
$T: \dom(T) \to X^*$ 
satisfies the properties:   
$\dom(T)$ is a linear subset of $X_d^*$ containing $\lin\{\delta_i^*\}$, 
$T(\sum_{i=1}^n c_i \delta_i^*)\to T(\sum_{i=1}^\infty c_i \delta_i^*)$ as $n\to\infty$ for every $\sum_{i=1}^\infty c_i \delta_i^*\in\dom(T)$, and there exists $\lambda\in(0,\infty)$ so that $\|T^*(F)\|\geq \lambda \|F\|$ for every $F\in\dom(T^*)$.

\item[{\rm (iii)}] The $X_d$-frames for $X$ are precisely the sequences $(T \delta_i^*)_{i=1}^\infty$,
where $T:X_d^*\to X^*$ is a bounded surjective operator.

\item[{\rm (iv)}] The Banach frames for $X$ wrt $X_d$ are precisely the sequences $(T \delta_i^*)_{i=1}^\infty$, where $T:X_d^*\to X^*$ is a bounded surjective operator which has a bounded right inverse defined from $X^*$ into $X_d^*$.

\item[{\rm (v)}] The sequences $\mg$ for which there exists an atomic decomposition  $(\mg, \mf)$ for $X$ wrt $X_d$ 
are precisely the sequences $(T \delta_i^*)_{i=1}^\infty$,
where $T:X_d^*\to X^*$ is a bounded surjective operator such that 

\begin{itemize}
\item[$({\mathcal P}_1):$]
 $T^*$ has a left inverse operator 
$L: \dom(L) \to X$ with $\dom(L)$ being a linear subset of $X_d$ 
 satisfying 
$ \dom(L)\supseteq \lin\{\delta_i\}\cup \ran(T^*)$ and 
$L(\sum_{i=1}^n T \delta_i^*(f) \delta_i)\to L(\sum_{i=1}^\infty T \delta_i^*(f) \delta_i)$ as $n\to\infty$ 
 for every $f\in X$. 
\end{itemize}

\item[{\rm (vi)}] The sequences $\mf$ for which there exists an atomic decomposition  $(\mg, \mf)$ for $X$ wrt $X_d$
are precisely the sequences $(T \delta_i)_{i=1}^\infty$,
where the operator  
$T: \dom(T) \to X$ 
satisfies the property
 \begin{itemize}
\item[$({\mathcal P}_2):$]
 $\dom(T)$ is a linear subset of $X_d$ 
 containing $\lin\{\delta_i\}$, 
$T(\sum_{i=1}^n c_i \delta_i)\to T(\sum_{i=1}^\infty c_i \delta_i)$ as $n\to\infty$ for every $\sum_{i=1}^\infty c_i \delta_i\in\dom(T)$ 
and $T$ has a bounded right inverse $U: X\to X_d$ with $\ran(U)\subseteq \dom(T)$,  
\end{itemize}
and such that $U^*$ is surjective. 

\item[{\rm (vii)}] 
The $X_d$-Riesz bases for $X$ are precisely the sequences $(T \delta_i)_{i=1}^\infty$, where $T:X_d\to X$ is a bounded bijective operator.
\end{itemize}

\end{thm}

\begin{proof}  (i)  follows easily from \cite[Cor. 3.3]{CCS}.

(ii) First assume that the operator $T: \dom(T)\to X^*$ satisfies the conditions listed in (ii) and consider the sequence $\mg$ given by $g_i=T \delta_i^*, i\in\mn$. 
Let $f\in X$ be such that $(g_i(f))\in X_d$ and let $F$ denote its corresponding element in $X^{**}$. 
Consider an arbitrary element $\sum_{i=1}^\infty c_i \delta_i^* \in \dom(T)$.  
Let $C$ denote the basis constant of the canonical basis of $X_d$. 
For every $n\in\mn$, $\sum_{i=1}^n c_i \delta_i^*\in \dom(T)$ and
\begin{eqnarray}\label{eqft}
 \left|FT\left(\sum_{i=1}^n c_i \delta_i^*\right) \right| &=& 
 \left|\sum_{i=1}^n c_i g_i(f) \right|  
 \leq    C \left\|\sum_{i=1}^n c_i \delta_i^* \right\|_{X_d^*}\, \left\|\sum_{i=1}^\infty g_i(f) \delta_i\right\|_{X_d}.
 \end{eqnarray}
 Using the assumptions on $T$ and taking limit  when $n\to\infty$, 
we obtain 
\begin{eqnarray*}
 \left|FT\left(\sum_{i=1}^\infty c_i \delta_i^*\right) \right|   
 \leq    C \left\|\sum_{i=1}^\infty c_i \delta_i^* \right\|_{X_d^*}\, \left\|\sum_{i=1}^\infty g_i(f) \delta_i\right\|_{X_d}.
 \end{eqnarray*}
Therefore, $FT$ is bounded on $\dom(T)$ and thus $F$ belongs to $\dom(T^*)$.  
Furthermore,  $(g_i(f))_{i=1}^\infty =(T^*(F)  (\delta_i^*))_{i=1}^\infty\in X_d$ 
and 
$$\|(g_i(f))_{i=1}^\infty\|_{X_d} = 
 \|T^*(F)\|\geq 
 \lambda \|f\|_X,$$
 which completes the proof.
 
 Conversely, assume that $\mg$ satisfies the lower $X_d$-frame condition.  Then the operator $T=T_\mg$ has the required properties.  
 
 (iii) follows easily from  \cite[Th. 3.9]{Sduals}.
 
(iv) Let $\mg$ be a Banach frame for $X$ w.r.t. $X_d$ and let $Q$ denote a Banach frame operator for $\mg$.  
By (iii),  
the operator $T_\mg:X_d^*\to X^*$ is bounded and surjective.
Furthermore, 
$Q^*$ is a bounded right inverse of $T_\mg$.  
Then the operator $T=T_\mg$ has the desired properties. 

Conversely, assume that  $T:X_d^*\to X^*$ is a bounded surjective operator which has a bounded right inverse 
$W:X^*\to X_d^*$. 
Consider the sequence $(g_i)_{i=1}^\infty =(T\delta_i^*)_{i=1}^\infty$. 
By (iii), $\mg$ is an $X_d$-frame for $X$. 
Furthermore, $W^*$ is a Banach frame operator for $\mg$.

(v) Let  $\mg$ be such that there exists an atomic decomposition $(\mg, \mf)$ for $X$ wrt $X_d$. 
By (iii), 
$T_\mg$ is bounded and surjective. 
By Lemma \ref{lc},  $\mf$ satisfies the lower $X_d^*$-frame condition. 
Consider the synthesis operator  
\begin{eqnarray}
& & T_\mf:\dom(T_\mf) \to X, \ \ T_\mf \seqgr[c]=\sum_{i=1}^\infty c_i f_i, \label{a1}\\ 
&\mbox{where }& \dom(T_\mf) =\left\{ \seqgr[c]\in X_d : \sum_{i=1}^\infty c_i f_i \mbox{ converges in } X\right\}.\label{a2}
\end{eqnarray}  
Clearly,  $\dom(T_\mf)$ is a linear set containing $ \lin\{\delta_i\}$. 
For every $f\in X$, one has $f=\sum g_i(f) f_i$,  implying that 
 $R(T^*_\mg)=R(U_\mg)\subseteq \dom (T_\mf)$ and that for every $f\in X$, 
 $T_\mf T_\mg^* (f)  
 =f$
and
$ T_\mf\left(\sum_{i=1}^n g_i(f) \delta_i\right) 
{\underset{n\to\infty}{\xrightarrow{\hspace*{.5cm}}}} \
\sum_{i=1}^\infty g_i(f) f_i
=  T_\mf \left(\sum_{i=1}^\infty g_i(f) \delta_i\right)$. 
 Take $T=T_\mg$ and $L=T_\mf$. 

Conversely, assume that the operator $T$ satisfies the conditions listed in (v) and consider the sequence $\mg=(T\delta_i^*)_{i=1}^\infty$. By (iii), $\mg$ is an $X_d$-frame for $X$. Define $f_i:=L\delta_i$, $i\in\mn$. For every $f\in X$ and $n\in\mn$,  
\begin{equation}\label{vconv}
\sum_{i=1}^n g_i(f) f_i = 
L \left(\sum_{i=1}^n g_i(f) \delta_i\right) \
{\underset{n\to\infty}{\xrightarrow{\hspace*{.5cm}}}} \
 L\left(\sum_{i=1}^\infty g_i(f) \delta_i\right) = 
 L T^*(f)=f,
 \end{equation}
which completes the proof  that  $(\mg, \mf)$ is an atomic decomposition for $X$ wrt $X_d$. 

\sloppy
(vi) Let $\mf$ be such that an atomic decomposition $(\mg, \mf)$ for $X$ wrt $X_d$ exists.  
  By Lemma \ref{lc},  $\mf$ satisfies the lower $X_d^*$-frame condition. 
Consider the operator $T_\mf$ given by (\ref{a1}) and (\ref{a2}). 
 Clearly, $\dom(T_\mf)\supseteq \lin\{\delta_i\}$ 
and for every $\seqgr[c]\in\dom(T_\mf)$, we have $T_\mf(\sum_{i=1}^n c_i \delta_i)
=\sum_{i=1}^n c_i f_i 
\to 
T_\mf(\sum_{i=1}^\infty c_i \delta_i)$ as $n\to\infty$.

 Since $\seqgr[g]$ is an $X_d$-frame for $X$, it follows that $U_\mg$ is bounded with bounded inverse on $\ran(U_\mg)$ which implies (see, e.g., \cite[Th. 4.15]{Rbook}) that $U_\mg^*$ is surjective. 
 Furthermore, 
 we have $\ran(U_\mg)\subseteq \dom(T_\mf)$ and  $U_\mg$ is a right inverse of $T_\mf$. 
Take $T=T_\mf$ and $U=U_\mg$.

Conversely, assume that $T$ satisfies the conditions listed in (vi) and consider the sequence 
$f_i=T\delta_i$, $i\in\mn$.
Define $g_i:= U^* \delta_i^*$, $i\in\mn$. By (iii), $\mg$ is an $X_d$-frame for $X$. For every $f\in X$,  we have that 
$(g_i(f))_{i=1}^\infty =Uf \in \ran(U)\subseteq \dom(T)$ and 
$$ \sum_{i=1}^n g_i(f) f_i = T\left(\sum_{i=1}^n g_i(f) \delta_i\right) 
\to 
T\left(\sum_{i=1}^\infty g_i(f) \delta_i\right) =
TUf = f,$$
which completes the proof.

(vii) follows easily from \cite[Prop. 3.4]{Srbasis}.
\end{proof}

Clearly, the characterizations in the above theorem concern the synthesis operators of the corresponding sequences,
i.e., 
one can write that  a sequence 
 is of type (i) (resp. (ii), ..., (vii)) if and only if its synthesis operator  
satisfies the properties 
of the operator $T$ stated 
in (i) (resp. (ii), ..., (vii)).  
Notice that the characterization of sequences satisfying the lower $X_d$-frame condition (Theorem \ref{allfr}(ii)) can be 
simplified if one uses an operator 
defined just on $\lin\{\delta_i^*\}$:

\begin{prop} \label{charlower2}
The sequences $G$ satisfying the lower $X_d$-frame condition 
are precisely the sequences $(T \delta_i^*)_{i=1}^\infty$,
where the 
operator  
$T: \lin\{\delta_i^*\} \to X^*$ 
has the property that there exists $\lambda\in(0,\infty)$ so that $\|T^*(F)\|\geq \lambda \|F\|$ for every $F\in\dom(T^*)$.
\end{prop}
\begin{proof} 
If $\mg$ satisfies the lower $X_d$-frame condition, the operator $T=T_\mg\mid_{\lin\{\delta_i^*\}}$ has the desired property. 
 
Conversely, assume that the densely defined operator $T:  \lin\{\delta_i^*\}\to X^*$ satisfies the property stated in the proposition and consider the sequence $\mg$ given by $g_i=T \delta_i^*, i\in\mn$. 
 Take $f$, $F$ and $C$ as in the proof of Theorem \ref{allfr}(ii). For every element $\sum_{i=1}^n c_i \delta_i^*\in \lin\{\delta_i^*\}$,
(\ref{eqft}) holds, 
implying that $FT$ is bounded on $\lin\{\delta_i^*\}$ and hence, $F$ belongs to $\dom(T^*)$. 
The conclusion for the lower inequality of 
$(g_i(f))_{i=1}^\infty$ follows as in Theorem \ref{allfr}(ii). 
 \end{proof}

To see that the above simplification of Theorem \ref{allfr}(ii)  is essential (even in the Hilbert space setting),  consider  for example the sequence $\mf=(ne_n)_{n=1}^\infty$ which satisfies the lower frame condition and  
$\dom(T_\mf)$ 
 is strictly larger then $ \lin\{e_i\}$. 
 In this manner, 
 Proposition \ref{charlower2} in the context of Hilbert spaces simplifies \cite[Prop. 4.6(e)]{BSA}.

\subsection*{Operators keeping the sequence type} 

For a given sequence $\mg$, here we consider conditions on an operator $V$ defined on $\lin \{g_i\}$ which preserve the sequence type. 
First, let us consider the case of $X_d$-Bessel sequences.
If $\mg$ is an $X_d$-Bessel sequence for $X$, then $(Vg_i)_{i=1}^\infty$ is an $X_d$-Bessel sequence for $X$ if and only if the operator $VT_\mg\vert_{\lin\{\delta_i^*\}}$ is bounded.
Note that the condition \lq\lq the operator $VT_\mg\vert_{\lin\{\delta_i^*\}}$ is bounded\rq\rq\ can not be relaxed to the condition
\lq\lq the operator $VT_\mg\vert_{\lin\{\delta_i^*\}}$ is bounded on the set $\{\delta_i^*\}_{i=1}^\infty$\rq\rq. Consider for example the Bessel sequence $\mg=(e_i)_{i=1}^\infty$ 
and the operator $V$ given by $Ve_i:=e_1$, $i\in\mn$.
In a similar way as in the $X_d$-Bessel sequence case, using Theorem \ref{allfr} one can list conditions on $V$ (more precisely, on the synthesis operator $T_{(Vg_i)}$) which are necessary and sufficient to preserve the type of the sequences discussed in Theorem \ref{allfr} and we will skip the listing.
In general, these conditions do not require $V$ to be bounded. 
Consider for example the Bessel sequence 
$(g_i)_{i=1}^\infty=(\frac{1}{i}e_i)_{i=1}^\infty$ 
(or the sequence 
$(g_i)_{i=1}^\infty=(ie_i)_{i=1}^\infty$ satisfying the lower frame condition) and the operator $V$ given by $Ve_i:=ie_i, i\in\mn$. 
However, it is not difficult to observe that operators which keep the $X_d$-Riesz basis property must be bounded:

\begin{prop}
Let $X_d$ be a $CB$-space and let $\mf$ be an $X_d$-Riesz basis for $X$. Let an operator $V: \lin\{f_i\} \to X$ be given. If $(Vf_i)_{i=1}^\infty$ is an $X_d$-Riesz basis for $X$, then $V$ is bounded.
\end{prop}

 Below we concentrate on bounded operators and determine the additional conditions which are necessary and sufficient to preserve the sequence type.
 We are motivated by \cite{Aldr}, where the author investigates constructions of frames using a given frame. 
This is of interest for applications where one aims at construction of frames with desired suitable properties. 
Note that Proposition \ref{p1}(ii) extends \cite[Theorem 1]{Aldr}.

\begin{prop} \label{p1}
Let $X$ be reflexive, $X_d$ be an $RCB$-space and $V$ be a bounded operator 
from $X^*$ into $X^*$  for {\rm (i)-(iv)} and  from $X$ into $X$ for {\rm (v)}.  
Then the following statements hold.
\begin{itemize}

\item[{\rm (i)}] Let $\mg$ be an $X_d$-Bessel sequence for $X$. Then 
$(Vg_i)_{i=1}^\infty$ is an $X_d$-Bessel sequence for $X$.

\item[{\rm (ii)}] Let $\mg$ be an $X_d$-frame for $X$. The sequence
 $(V g_i)_{i=1}^\infty$ is an $X_d$-frame for $X$ 
if and only if $V$ is
surjective.

\item[{\rm (iii)}]  Let $\mg$ be a Banach frame for $X$ wrt $X_d$. The sequence 
 $(V g_i)_{i=1}^\infty$ is a Banach frame for $X$ wrt $X_d$ 
if and only if $V$ is surjective and  has a bounded right inverse $W:X^*\to X^*$.

\item[{\rm (iv)}] Let $\mg$ be such that there exists an atomic decomposition $(\mg, \mf)$ for $X$ wrt $X_d$.
The sequence $(V g_i)_{i=1}^\infty$ is the first one in some atomic decomposition pair if and only if $V$ is surjective and $U_\mg V^*$ has a left inverse $L:\dom(L)(\subseteq X_d)\to X$ with $\dom(L)$ being a linear subset of $X_d$ satisfying 
$ \dom(L)\supseteq \lin\{\delta_i\}\cup \ran(U_\mg V^*)$ and 
$LV(\sum_{i=1}^n g_i(f) \delta_i)\to LV(\sum_{i=1}^\infty g_i(f) \delta_i)$ as $n\to\infty$ 
 for every $f\in X$. 

\item[{\rm (v)}] 
Let $\mf$ be an $X_d$-Riesz basis for $X$. The sequence $(V f_i)_{i=1}^\infty$ is an $X_d$-Riesz basis for $X$ if and only if $V$ is bijective. 
\end{itemize}
\end{prop}

\begin{proof}  We sketch the proofs.

(i) If $B_\mg$ denotes an $X_d$-Bessel bound for $\mg$, then for every $f\in X$ one has $((Vg_i)(f))_{i=1}^\infty=(g_i(V^*f))_{i=1}^\infty\in X_d$ and
$\| ((Vg_i) (f))_{i=1}^\infty\|_{X_d}\leq  B_\mg \, \|V\| \cdot\|f\|$.

 (ii) 
By (i),  $(V g_i)_{i=1}^\infty$ is an $X_d$-Bessel sequence for $X$. 
Furthermore,  $(V g_i)_{i=1}^\infty $ satisfies the lower $X_d$-frame condition
 if and only if there exists $\lambda>0$ so that $\|V^*f\|\geq \lambda \|f\|$, $\forall f\in X$, which is known to be equivalent to $V$ being surjective. 

(iii) Let $Q$ denote a Banach frame operator for $\mg$.

First assume that $V$ is surjective and  has a bounded right inverse $W:X^*\to X^*$. 
By (ii), $\mg$ is an $X_d$-frame for $X$. Furthermore, $W^*Q$ is a Banach frame operator for $(Vg_i)_{i=1}^\infty$.

Conversely, assume that $(Vg_i)_{i=1}^\infty$ is a Banach frame for $X$ wrt $X_d$ and let 
$Q_1$ denote a Banach frame operator for $(Vg_i)_{i=1}^\infty$. By (ii), $V$ is surjective. Furthermore, 
$T_\mg (Q_1)^*$  
 is a bounded right inverse of $V$.

(iv) and (v)  follow easily using Theorem \ref{allfr}. 
\end{proof}

\begin{rem} For a given $X_d$-frame (resp. $X_d$-Bessel sequence, Banach frame, atomic decomposition, sequence satisfying the lower $X_d$-frame condition) $\mg$ for
$X$, not all the $X_d$-frames (resp. $X_d$-Bessel sequences, Banach frames, atomic decompositions, sequences satisfying the lower $X_d$-frame condition) for $X$ can be obtained 
in the way $(Vg_i)$ using an operator $V$. 
Consider for example the frame $\mg=(e_1, e_1, e_2, e_3, e_4, \ldots)$ and the frame $(e_1, e_2, e_3, e_4, e_5, \ldots)$ for $\h$, which can not be written as $(Vg_i)$ for any operator $V$.
In contrary, when $X_d$ is a $CB$-space
and an $X_d$-Riesz basis $\mf$ for $X$ is given, then every $X_d$-Riesz basis $\mw=(w_i)$ for $X$ can be 
written in the way $(Vf_i)$ using the operator $V=T_\mw T_\mf^{-1}$ which is a bounded bijection of $X$ onto $X$. 
\end{rem}

\section{Series expansions in a Banach space and its dual space via $X_d$-frames}\label{s4}

While a Hilbert frame $\seqgr[g]$ for $\h$ always has a sequence $\seqgr[f]$ satisfying (\ref{fexp2}), this is not 
always  
the case with an $X_d$-frame.
The sequence $\seqgr[g]=(e_i+e_{i+1})_{i=1}^\infty$  is an $X_d$-frame (even  a Banach frame) for $\h$ with respect to an 
 appropriate sequence space $X_d$ and there is no sequence $\seqgr[f]\in \h^\mn$ such that $f=\sum_{i=1}^\infty \langle f,g_i\rangle  f_i$ holds for all $f\in \h$ \cite{CCS}. 
Theorem \ref{allfr}(v) gives a characterization of the $X_d$-frames $\mg$ for $X$ for which there exists $\mf$ satisfying  
\begin{equation} \label{fx}
f=\sum_{i=1}^\infty g_i(f) f_i, \ \forall f\in X.
\end{equation}
Let $\mg$ be an $X_d$-frame for $X$. By  \cite{CCS}, when $X_d$ is an $RCB$-space 
 and $\ran(U_\mg)$ is complemented in $X_d$, then there exists $\mf$ which is an $X_d^*$-frame for $X^*$ and satisfies (\ref{fx}); such $\mf$ is called a {\it dual $X_d^*$-frame of $\mg$}. 
  In some cases, there might exist a sequence $\mf$ which is not  an $X_d^*$-frame for $X^*$ but yet satisfies  (\ref{fx}), see Example \ref{adn}. 
Every $\mf$ satisfying (\ref{fx}) will be called a {\it synthesis-pseudo-dual of $\mg$}, or in short, {\it s-pseudo-dual of $\mg$}, as an analogue to the Hilbert space case. Theorem \ref{pr51} below gives a characterization of all
the s-pseudo-duals of $\mg$, if such ones exist, and a characterization of those ones among them, which are $X_d^*$-frames for $X^*$.

\begin{ex} \label{adn} {\rm \cite[Ex. 5.1]{Sduals}} 
Let $X=\ell^p$ ($1<p<\infty$), $\seqgr[\xi]$ be the canonical basis of
$X$ and $\seqgr[E]$ be the coefficient functionals associated to
$\seqgr[\xi]$. Then the sequence $\mg = \left(\frac{1}{2} E_1, E_2,
\frac{1}{2^2} E_1, E_3, \frac{1}{2^3} E_1, E_4, \ldots \right)$ is a
$p$-frame for $X$,
the sequence $\mf = (\xi_1, \xi_2, \xi_1, \xi_3,
\xi_1, \xi_4, \ldots)$ satisfies the lower $q$-frame condition, but does not satisfy the upper one, and
\begin{equation*} \label{ex2}
f=\sum_{i=1}^\infty g_i(f)f_i, \ \forall f\in X, \ \ \mbox{and} \ \ g=\sum_{i=1}^\infty
g(f_i)g_i, \ \forall g\in X^*.
\end{equation*}
\end{ex}

\begin{thm} \label{pr51}
Let $X_d$ be an $RCB$-space and let $\mg$ be an $X_d$-frame for $X$. Then the following holds. 
 \begin{itemize}
 \item[{\rm (a)}] If $\mg$ has a dual $X_d^*$-frame, 
then all the dual $X_d^*$-frames
  of $\mg$ are precisely the sequences $(L \delta_i)_{i=1}^\infty$,
where $L:X_d\to X$ is a bounded linear extension of $U_\mg^{-1}$.

\item[{\rm (b)}] If $\mg$ has an s-pseudo-dual, 
 then all the s-pseudo-duals of $\mg$ are precisely the sequences  $(L\delta_i)_{i=1}^\infty$ where  the operator 
 $L:\dom(L) \to X$ is such that

 \begin{itemize}
\item[$({\mathcal P}_3):$]  $\dom(L)$ is a linear subset of $X_d$ containing $\lin \{\delta_i\}$ and $\range(U_\mg) $, $L(\sum_{i=1}^n g_i(f) \delta_i)\to L(\sum_{i=1}^\infty g_i(f) \delta_i)$ as $n\to\infty$ for every $f\in X$,
and 
$L$ is a left inverse of $U_\mg$. 
\end{itemize}

\end{itemize}
\end{thm}
\begin{proof} The assumptions imply that $X$ is isomorphic to $\range(U_\mg)$ and $\range(U_\mg)$ is a closed subspace of $X_d$ \cite{CS,Sduals}, and thus $X$ is also reflexive.

(a) Let $\mg$ have a dual $X_d^*$-frame. By  \cite[Prop. 3.4]{CCS}, $U_\mg^{-1}$ can be extended to a bounded operator on $X_d$. 
If $L$ is a bounded linear extension of $U_\mg^{-1}$ on $X_d$, it is proved 
in \cite[Prop. 3.4]{CCS} that $(L \delta_i)_{i=1}^\infty$ is a dual $X_d^*$-frame of $\mg$. Conversely, let $\mf$ be a dual $X_d^*$-frame of $\mg$.
By Theorem \ref{allfr}(iii), 
 the sequence $\mf$ has the form $(T \delta_i)_{i=1}^\infty$ for a bounded surjective operator $T:X_d\to X$. Furthermore, (\ref{fx}) implies that
$T$ is an extension of $U_\mg^{-1}$. Take $L=T$. 

(b) Assume that $\mg$ has an s-pseudo-dual $\mf$. 
In a similar way as in the proof of  Theorem \ref{allfr}(v), consider the operator $T_\mf$ given by  (\ref{a1}) and (\ref{a2}) and observe that it satisfies the properties of the wanted operator $L$, so take $L=T_\mf$.

Conversely, consider a sequence $(L\delta_i)_{i=1}^\infty$ with $L$ satisfying the properties described in (b). Then for every $f\in\h$, (\ref{vconv}) holds, 
and hence,  $(L\delta_i)_{i=1}^\infty$ is an s-pseudo-dual of $\mg$. 
 \end{proof}

Thus, when $X_d$ is an $RCB$-space and $\mg$ is an $X_d$-frame for $X$, then an s-pseudo-dual $\mf$ of $\mg$ is a dual $X_d^*$-frame of $\mg$ if and only if  $T_\mf$ is bounded on $\dom(T_\mf)$. 
For conditions equivalent to the existence of dual $X_d^*$-frames of $X_d$-frames see \cite{CCS, Sb}.

\vspace{.1in}
There exist series expansions in the form (\ref{fx}) where $\mg$ is an $X_d$-Bessel sequence and not an $X_d$-frame. 
As a trivial example, consider the Bessel sequence $\mg=(\frac{1}{i}e_i)_{i=1}^\infty$ in $\h$ and the expansions $f=\sum_{i=1}^\infty \langle f, \frac{1}{i}e_i\rangle i e_i$ valid for all $ f\in\h$. 
Notice that some characterizations from Theorems \ref{allfr} and \ref{pr51} related to $X_d$-frames  can  easily be extended to  $X_d$-Bessel sequences:

\begin{prop} \label{pb}  Let $X_d$ be an $RCB$-space and $X$ be reflexive. Then the following statements hold. 
 \begin{itemize}
 \item[{\rm (a)}] The $X_d$-Bessel sequences $\mg$ for $X$, for which there exists $\mf$ satisfying (\ref{fx}), are precisely the sequences 
$(T \delta_i^*)_{i=1}^\infty$, 
where $T:X_d^*\to X^*$ is a bounded operator such that 
 $({\mathcal P_1})$ holds.

 \item[{\rm (b)}] Let $\mg$ be an $X_d$-Bessel sequence for $X$. If there exists $\mf$ satisfying (\ref{fx}), 
 then all possibilities for $\mf$ 
 are precisely the sequences  $(L\delta_i)_{i=1}^\infty$ where the operator 
$L:\dom(T)\to X$ is such that $({\mathcal P}_3)$ holds.

   \item[{\rm (c)}] The sequences $\mf$ for which there exists an $X_d$-Bessel sequence $\mg$ for $X$ satisfying (\ref{fx})  are precisely the sequences $(T \delta_i)_{i=1}^\infty$,
where the operator 
$T:\dom(T)\to X$ satisfies $({\mathcal P}_2)$.

\end{itemize}
\end{prop}

Let $X_d$ be an $RCB$-space. Note that there is no need to characterize dual $X_d^*$-Bessel sequences of $X_d$-Bessel sequences, because if $\mg$ is an $X_d$-Bessel sequence for $X$ and $\mf$  is an $X_d^*$-Bessel for $X^*$ satisfying (\ref{fx}), then $\mg$ is an $X_d$-frame for $X$ 
and $\mf$ is an $X_d^*$-frame for $X^*$
\cite{Sduals}.

\subsection*{Connection between expansions in $X$ and $X^*$}

Let $X_d$ be an $RCB$-space 
and let $\mg$ be an $X_d$-Bessel sequence for $X$.
If $\mf$ is an $X_d^*$-Bessel sequence for $X^*$, then $\mf$ satisfies (\ref{fx}) if and only if it satisfies
\begin{equation} \label{fg}
g=\sum_{i=1}^\infty g(f_i) g_i, \ \forall g\in X^*,
\end{equation}
see \cite[Lemma 4.3]{Sduals}. 
If $\mf$ is not an $X_d^*$-Bessel sequence for $X^*$, such equivalence does not hold in general. 
There exist cases where both representations (\ref{fx})  and (\ref{fg}) hold 
 (see Example \ref{adn}) and
there exist cases where (\ref{fx}) holds, but (\ref{fg}) does not hold (see Example \ref{eximp}).  
In Theorem \ref{1s2} below we determine a subset of $X^*$ where the representation $g=\sum_{i=1}^\infty g(f_i)g_i$ holds based on validity of (\ref{fx}). 
For investigation of a converse situation, namely, conclusions for $f=\sum g_i(f)f_i$ 
assuming validity of $g=\sum_{i=1}^\infty g(f_i)g_i, \ \forall g\in\dom(U_\mf)$, see \cite[Theor. 5.1]{Slower} switching the roles of $(f_i)_{i=1}^\infty$ and $(g_i)_{i=1}^\infty$,  and  switching the roles of $X$ and $X^*$.

\begin{thm} \label{1s2}
Let $X_d$ be an $RCB$-space, $\mg$ be an $X_d$-Bessel sequence for $X$, and $\mf$ satisfy (\ref{fx}). Then 
$g=\sum_{i=1}^\infty g(f_i) g_i$  for every $g\in D(T_\mf^*)$.
\end{thm}
\begin{proof}  By Lemma \ref{lc}, the sequence $\mf$ satisfies the lower $X_d^*$-frame condition. 
Consider the operator $T_\mf$ given by (\ref{a1}) and (\ref{a2}).  Clearly, $\mf=(T_\mf \delta_i)_{i=1}^\infty$ and by (\ref{fx}), $T_\mf$ is a left inverse of $U_\mg$. 
Furthermore, for every $g\in D(T_\mf^*)$ we have
$\sum_{i=1}^n g(f_i) g_i = T_\mg(\sum_{i=1}^n gT_\mf(\delta_i) \,\delta_i^*)
\underset{n\to\infty}{\longrightarrow}
T_\mg T_\mf^* (g) =
U_\mg^*T_\mf^* (g)=g. $
\end{proof}

\begin{rem} \label{setopt} The set $D(T_\mf^*)$ in the conclusion of the above proposition is the optimal one.  
The sequences in Example \ref{eximp} fulfill the assumptions of Theorem \ref{1s2}, 
$D(T_\mf^*)\neq \h$, 
and 
$$g=\sum_{i=1}^\infty \langle g, f_i\rangle  g_i \ \mbox{ if and only if } \ g\in D(T_\mf^*).$$ 
This shows that Theorem \ref{1s2} can not be improved in the sense that 
the set of validity of $g=\sum_{i=1}^\infty g(f_i) g_i$ can not be extended in general. 
However, there exist certain cases where the assumptions of Theorem \ref{1s2} hold,  $D(T_\mf^*)\neq X^*$, and the representation 
$g=\sum_{i=1}^\infty \langle g, f_i\rangle  g_i $ hold for every $g\in X^*$ - see Example \ref{adn}, where we have that  (\ref{fx}) is satisfied, $E_1\in X^*\setminus  D(T_\mf^*)$ and (\ref{fg}) holds. 
\end{rem}

\begin{rem} \label{setopt2} Given sequence $\mf$, consider the operator $T_\mf$, given by (\ref{a1}) and (\ref{a2}), and the operator $U_\mf:\dom(U_\mf) \to X_d^*$ determined by $U_\mf g=(g(f_i))_{i=1}^\infty$ for $g\in \dom(U_\mf) =\left\{ g\in X^*\  : \ (g(f_i))_{i=1}^\infty \in X_d^*\right\}$. 
 Notice that when $X_d=\ell^p$, $p\in (1,\infty)$, the set $\dom(T_\mf^*)$ in the above theorem coincides with $\dom(U_\mf)$.
 Actually, for any sequence $\mf$ and $X_d=\ell^p$, $p\in (1,\infty)$, one has that   $T_\mf^*=U_\mf$ - 
when $X$ is a Hilbert space and $X_d=\ell^2$, this is proved in \cite[Prop. 3.39(i)]{BSA} and the proof can be easily extended to the case when $X$ is a Banach space and $X_d=\ell^p$, $p\in (1,\infty)$, 
using \cite[Ex. 34.2]{Heuser}. 
\end{rem}

\section{The case of Hilbert frames and their dual sequences not necessarily being frames} \label{s5}

Here we apply some results from the previous sections to the Hilbert space setting and furthermore deepen the investigation, focusing on frames and their dual sequences (not necessarily frames) being important for reconstructions in applications. 
As we show via an example, interchanging a frame and a non-frame sequence in series representations in the form  (\ref{fexp2}) and (\ref{fexp1}) is not always possible, which leads us to further investigation of each of the two types of representations.

\subsection*{Connections between the representations (\ref{fexp2}) and (\ref{fexp1}) } 
\label{conn12}

Let $\seqgr[g]$ be a Bessel sequence in $\h$. It is well known that a Bessel sequence $\seqgr[f]$ in $\h$ satisfies 
(\ref{fexp2})  if and only if it satisfies (\ref{fexp1}) and in this case both  $\seqgr[g]$ and  $\seqgr[f]$ are necessarily frames for $\h$ (see, e.g., \cite[Lemma 6.3.2]{Cbook}).  
When $\seqgr[g]$ is a frame for $\h$ and $\seqgr[f]$ is not a Bessel sequence in $\h$, there are still certain cases where both (\ref{fexp2}) and (\ref{fexp1}) may hold (see Example \ref{adn} with $p=2$), but in general it is not possible to state an equivalence of (\ref{fexp2}) and (\ref{fexp1}):

\begin{ex} \label{eximp}
Consider the frame 
$\seqgr[g]=(e_1, e_1, e_1, e_2, e_2, e_2, e_3, e_3, e_3, \ldots)$ for $\h$ and the 
sequence 
$\seqgr[f]=(e_1, e_1, -e_1, e_2, e_1, -e_1, e_3, e_1, -e_1,\ldots)$ which satisfies the lower frame condition but is not Bessel in $\h$. Then 
(\ref{fexp2}) holds, 
but (\ref{fexp1}) does not hold (for example, $g=e_1$ breaks (\ref{fexp1})). 
Furthermore, the representation $g=\sum_{i=1}^\infty \langle  g,f_i \rangle g_i$ holds if and only if $g$ belongs to the closed linear span of $\{e_i\}_{i=2}^\infty$ if and only if 
 $g\in \dom(T_\mf^*)$ if and only if $g\in\dom(U_\mf)$.
 \end{ex}

The above example shows that \cite[Prop. 4.10]{LO} is not correct. 
 In  \cite[Prop. 4.10]{LO} it is claimed that when a frame  $\seqgr[g]$ for $\h$ and a sequence  $\seqgr[f]$ in $\h$ satisfy $\langle f,g\rangle = \sum_{i=1}^\infty \langle  f,g_i \rangle \langle  f_i, g\rangle$ for all $f,g\in\h$, then norm-convergence of (\ref{fexp2}) is equivalent to norm-convergence of  (\ref{fexp1}). However, this is not always the case as one can see in Example \ref{eximp}. 
  Note that if an unconditional convergence is considered, then there is an equivalence (even in cases where both $\seqgr[g]$ and $\seqgr[f]$ are not necessarily frames), more precisely, for any two sequences 
  $\seqgr[g]$ and $\seqgr[f]$, 
  (\ref{fexp2}) with unconditional convergence is equivalent to (\ref{fexp1}) with unconditional convergence \cite[Lemma 3.1]{uncconv2011}.

Assuming  validity of (\ref{fexp2}), the following statement determines the optimal subset of $\h$ where the representation $g=\sum_{i=1}^\infty \langle g, f_i\rangle g_i$ holds in general (though in particular cases it may hold for the entire space $\h$). It is a consequence of Theorem \ref{1s2} and the fact that $\dom(T_\mf^*)=\dom(U_\mf)$ for any sequence $\mf=(f_i)_{i=1}^\infty$ with elements from $\h$ \cite{BSA}.

\begin{cor}\label{1s2cor}
Let  $\seqgr[g]$ be a Bessel sequence in $\h$ and let $\mf=\seqgr[f]$ satisfy (\ref{fexp2}). Then 
 $g=\sum_{i=1}^\infty \langle  g,f_i \rangle g_i$   for every 
 $g\in\dom(U_\mf)$. 
  By Example \ref{eximp}, the set $\dom(U_\mf)$ is the optimal one in the sense that a general statement for a larger subset of $\h$ can not be claimed.
\end{cor}

\subsection*{The lower frame condition and its connection to expansions}

When  $\seqgr[g]$ is a frame for $\h$ and $\seqgr[f]$ satisfies (\ref{fexp2}) or (\ref{fexp1}), then $\seqgr[f]$ satisfies the lower frame condition for $\h$ \cite{CCLL,LO}.
The converse situation 
has also been of interest: given a sequence $\seqgr[f]$  satisfying the lower frame condition, 
what can be said about the existence of series expansions in the form (\ref{fexp2}) or (\ref{fexp1})?

\subsubsection*{Expansions in the form (\ref{fexp1})}

It is proved in \cite[Prop. 3.4]{CCLL} that {\it a sequence $\mf=(f_i)_{i=1}^\infty\in\h^\mn$ satisfies the lower frame condition if and only if there exists a Bessel sequence $(g_i)_{i=1}^\infty\in\h^\mn$ such that}
\begin{equation}\label{duf}
g=\sum_{i=1}^\infty \langle g, f_i\rangle  g_i, \ \forall g\in \dom(U_\mf).
\end{equation}  
In \cite{CCLL} and \cite{Slower} one can find examples which show that in general the representation $g=\sum_{i=1}^\infty \langle g, f_i\rangle  g_i$ does not need to be extendible 
to all $g\in\h$. 
Example \ref{eximp}  presents a case where the representation $g=\sum_{i=1}^\infty \langle g, f_i\rangle  g_i$ holds if and only if $g\in \dom(U_\mf)$.

\subsubsection*{Expansions in the form (\ref{fexp2})}

  Throughout this section, let  $\mf=(f_i)_{i=1}^\infty$ be a sequence with elements from $\h$ and consider its (densely defined) synthesis operator $T_\mf$.   
The consideration here 
 is motivated by the following result from \cite{Cinv}:
 
\begin{thm}\label{oleinv} {\rm \cite[Theor. 4.1]{Cinv}} Let $T_\mf$ be closed and surjective. Then $(f_i)_{i=1}^\infty$ satisfies the lower frame condition and there exists a Bessel sequence $(g_i)_{i=1}^\infty\in\h^\mn$ so that  (\ref{fexp2}) holds.

\end{thm}
Notice that when the densely defined  operator $T_\mf$ is closed, then it is surjective if and only if there is a positive number $a$ so that $\|T_\mf^* u\|\geq a\|u\|$ for every $u\in\dom(T_\mf^* )$ (see, e.g., \cite[Theor. 2.20]{Brezis}) which happens if and only if
 $\|U_\mf u\|\geq a\|u\|$ for every $u\in\dom(U_\mf)$, because  $T_\mf^*=U_\mf$ by \cite{BSA}. Thus, 
 Theorem \ref{oleinv} actually concerns precisely the  sequences $\mf$ which satisfy the lower frame condition and whose synthesis operator $T_\mf$ is closed. Here we consider sequences which satisfy the lower frame condition, but whose synthesis operator is not necessarily closed.

\begin{thm}  \label{lf} 
Let $(f_i)_{i=1}^\infty$ satisfy the lower frame condition and let $U_\mf$ be densely defined. Then
 there exists a Bessel sequence $(g_i)_{i=1}^\infty$ in $\h$ so that for $f\in\h$ one has that
\begin{equation*} \label{f3}
f=\sum_{i=1}^\infty \langle f, g_i\rangle  f_i \ \  
\mbox{if and only if} 
\ \sum_{i=1}^\infty \langle f, g_i\rangle  f_i \ \mbox{converges}.
\end{equation*}
\end{thm} 
\begin{proof}  
The statement can be derived as a consequence of \cite[Prop. 3.4]{CCLL} and \cite[Cor. 5.1]{Slower}. For the sake of completeness, we include a proof here. 

Since  $\ran(U_\mf)$  is a closed subspace of $\ell^2$ \cite[Lemma 3.1]{CCLL} and $U_\mf^{-1}$ (defined on $\ran(U_\mf)$)  is bounded, there exists a bounded extension $V:\ell^2\to\h$ of $U_\mf^{-1}$. Define $g_i=V\delta_i, i\in\mn$. Clearly, $\mg$ is a Bessel sequence in $\h$. 
Let $f\in\h$ be such that $\sum_{i=1}^\infty \langle f, g_i\rangle  f_i $ converges.
Using Theorem \ref{allfr}(ii) 
 and the fact that $T_\mf\subseteq U_\mf^* $ \cite[Prop. 4.6]{Cinv}, 
we have
$$ \sum_{i=1}^n \langle f,g_i\rangle f_i 
\to  T_\mf \left(\sum_{i=1}^\infty \langle f,g_i\rangle \delta_i\right) = 
U^*_\mf V^*(f)=f. $$
\end{proof}

\begin{rem}
Notice that Theorem \ref{lf} extends Theorem \ref{oleinv} in the following sense: {\it $(f_i)_{i=1}^\infty$ satisfies the assumptions of Theorem \ref{oleinv} if and only if it satisfies the assumptions of Theorem \ref{lf} and $T_\mf$ is closed.} Indeed, having in mind the arguments in the paragraph before Theorem \ref{lf}, 
it remains only  to show that 
Theorem \ref{oleinv} implies the density of $\dom(U_\mf)$ in $\h$. Under the assumptions of Theorem \ref{oleinv}, \cite[Theorem 13.12]{Rbook} implies that $T_\mf^*$ is densely defined, which implies that $U_\mf$ is densely defined, because  $T_\mf^*=U_\mf$ by \cite{BSA}. 

For a case where Theorem \ref{lf} applies and Theorem \ref{oleinv} does not apply, see Example \ref{ex1} below. For the purpose 
 of the example, first observe that  the sufficient condition for closedness of $T_\mf$ in \cite[Corol. 4.7]{Cinv} is also necessary:
\end{rem}

\begin{lem}\label{lemc1} 
The synthesis operator $T_\mf$ is closed if and only if $U_\mf$ is densely defined and $\dom(T_\mf)=\dom(U_\mf^*)$. 
\end{lem}
\bp
Assume that $T_\mf$ is closed. Then $\dom(T_\mf^*)$ is dense and $T_\mf=T_F^{**}$ (see, e.g., \cite[Theor. 13.12]{Rbook}). By \cite{BSA}, $U_\mf=T_\mf^*$, which now implies that  $U_\mf$ is densely defined and $U_\mf^*=T_\mf$. 
The other direction of the statement is given in \cite[Corol. 4.7]{Cinv}.
\ep

\vspace{.1in}

Now we consider an example where Theorem \ref{lf} applies with validity of (\ref{fexp2}), 
and where Theorem \ref{oleinv} does not apply.  
We use a small modification of the example from \cite[p. 413]{Cinv}.

\begin{ex} \label{ex1} Consider the sequence
$$\seqgr[f]=(e_1,  2(e_2-e_1), 2{\rm \small x}2e_2, 3(e_3-e_2), 3{\rm \small x}2e_3, 4(e_4-e_3),4{\rm \small x}2e_4, 5(e_5-e_4), \ldots).$$ Clearly, $\seqgr[f]$ satisfies the lower frame condition, but does not satisfy the upper one. Since $\sum_{i=1}^\infty  |\langle f_j, f_i\rangle |^2 <\infty$ for every $j\in\mn$,  $U_\mf$ is densely defined by \cite[Prop. 4.5]{Cinv}.  
Thus, $\seqgr[f]$ satisfies the assumptions of Theorem \ref{lf}.
Using the Bessel sequence $\seqgr[g]=(e_1, 0, \frac{1}{2{\rm \small x}2}e_2, 0, \frac{1}{3{\rm \small x}2}e_3, 0,\frac{1}{4{\rm \small x}2}e_4, 0,\ldots)$,  
 (\ref{fexp2}) and (\ref{fexp1}) hold.

Consider the sequence 
$\seqgr[c]=(1, 1/2, -1/2, -1/3, 1/3, 1/4, -1/4, -1/5, \ldots)$  which belongs to $\ell^2.$ 
Since the partial sums of the series 
$\sum_{i=1}^{\infty} c_i f_i$ are $e_1, e_2, -e_2, -e_3, e_3, e_4, -e_4, -e_5, \ldots$, 
it follows that the series $\sum_{i=1}^\infty c_i f_i$ does not converge and thus  $\seqgr[c]\notin \dom(T_\mf)$. 
For every $f\in\h$, the partial sums of the series
$\sum_{i=1}^{\infty}  \overline{c_i} \langle f,  f_i\rangle$ are
 $\langle f, e_1\rangle, \langle f, e_2\rangle, - \langle f, e_2\rangle, - \langle f,e_3\rangle, \langle f,e_3\rangle, \langle f,e_4\rangle, - \langle f,e_4\rangle, - \langle f,e_5\rangle, \ldots$, 
which implies that  $\sum_{i=1}^\infty\overline{c_i} \langle f,  f_i\rangle $ converges and equals  
$0 $; therefore, $\seqgr[c]\in \dom(U_\mf^*)$. 
Now Lemma \ref{lemc1} implies that $T_\mf$ is not closed and thus $\seqgr[f]$ does not satisfy the assumptions of Theorem \ref{oleinv}. 
\end{ex}

\subsection*{Characterization of all the s-pseudo-duals of a frame} \label{charactsduals}

As Example \ref{eximp} shows,  
in general (\ref{fexp2}) and (\ref{fexp1}) are not equivalent, which motivates 
their separate investigation. 
Given a frame $\seqgr[g]$ for $\h$, 
a characterization of the sequences $\seqgr[f]$ (not necessarily frames for $\h$) which satisfy (\ref{fexp2}) can be obtained as a consequence of Theorem  \ref{pr51}:

\begin{cor} \label{sdual}
Let $\mg=\seqgr[g]$ be a frame for $\h$. Then all the s-pseudo-duals of 
$\seqgr[g]$ are precisely the sequences  $(L\delta_i)_{i=1}^\infty$ where the operator $L:\dom(L) \to \h$ satisfies $({\mathcal P}_3)$ with $X_d=\ell^2$ and $X=\h$. 
\end{cor}

Recall that the dual frames of a frame $\mg=\seqgr[g]$ for $\h$ are characterized as the sequences $(L\delta_i)_{i=1}^\infty$ where the operator $L: \ell^2\to \h$ is 
 a bounded linear extension of $U_\mg^{-1}$ \cite[Lemma 6.3.5]{Cbook}. 
 Now Corollary \ref{sdual} clarifies the difference between a dual frame and an $s$-pseudo-dual which is not a frame, in terms of operator-properties - for s-pseudo-duals, the boundedness of $L$ is relaxed to the convergence-properties in $({\mathcal P}_3)$.

Notice that by \cite{LO}, given a frame $\mg$ for $\h$, the linear extensions $L$ of $U_\mg^{-1}$ can be characterized as the operators
\begin{equation}\label{w2}
L=S_\mg^{-1}T_\mg + W - WU_\mg S_\mg^{-1}T_\mg,
\end{equation} where  $W:\dom(W) (\subseteq \ell^2)\to \h$ is a linear operator with $\dom(W)\supseteq \ran(U_\mg)$. 
Using this characterization 
and Corollary \ref{sdual}, 
we can obtain s-pseudo-duals as follows:

\begin{cor} \label{sdual2}
Let $\mg=\seqgr[g]$ be a frame for $\h$. Let $W$ be a liner operator whose domain contains $\ran(U_\mg) \cup \lin\{\delta_i\}_{i=1}^\infty$
and such that

$\bullet$ $ W(\sum_{i=1}^n \langle f, g_i\rangle \delta_i) \to W(\sum_{i=1}^\infty \langle f, g_i\rangle \delta_i)$ as $n\to \infty$, $\forall f\in\h$,

$\bullet$ $W $ is bounded on $\ran(U_\mg) $.

\noindent Then the sequence $(L\delta_i)_{i=1}^\infty$, where $L$ is given by (\ref{w2}),  is an s-pseudo-dual of $\mg$. 
\end{cor}
\bp 
 Clearly,  
 $\dom(L)$ contains  $\ran(U_\mg)\cup\lin\{\delta_i\}_{i=1}^\infty$, 
and by \cite{LO}, 
$LU_\mg$ is the identity operator on $\h$. 
For $f\in\h$, using the assumptions on $W$, it follows that $(W - WU_\mg S_\mg^{-1}T_\mg) (\sum_{i=1}^n \langle f, g_i\rangle \delta_i)
\to (W - WU_\mg S_\mg^{-1}T_\mg)(\sum_{i=1}^\infty \langle f, g_i\rangle \delta_i)$ as $n\to\infty$ and since $S_\mg^{-1}T_\mg$ is continues, it leads to  $L(\sum_{i=1}^n \langle f, g_i\rangle \delta_i)
\to L(\sum_{i=1}^\infty \langle f, g_i\rangle \delta_i)$ as $n\to\infty$. 
Now apply Corollary \ref{sdual}. 
\ep

\vspace{.1in} As an illustration of Corollary \ref{sdual2}, consider the example below. 
Notice that the boundedness of $W$ on $\ran(U_\mg)$ is not equivalent to boundedness of $L$ on $\dom(L)$  (contrary to what is stated in \cite[p. 296]{LO}) and thus Corollary \ref{sdual2} does not limit the construction to s-pseudo-duals which are dual frames.

\begin{ex}
Consider the frame $\mg=(e_1, e_1, e_1, e_2, e_2, e_2, e_3, e_3, e_3, \ldots)$ for $\h$ and the operator $W$ defined by 

$\bullet$ $W\delta_{3k-2} = e_k, W\delta_{3k-1} = e_1, W\delta_{3k} =- e_1, k\in\mn$, and by lineariry on $\lin\{\delta_i\}_{i=1}^\infty$;

$\bullet$ $W(U_\mg f)=f$ for $U_\mg f\notin \lin\{\delta_i\}_{i=1}^\infty$ (possible because of the injectivity of $U_\mg$);
 for $U_\mg f\in \lin\{\delta_i\}_{i=1}^\infty$ - observe that calculating $W(U_\mg f)$ using the definition of $W$ on $\lin\{\delta_i\}_{i=1}^\infty$ gives $W(U_\mg f)=f$, so then for every $U_\mg f\in\ran(U_\mg)$ one has $W(U_\mg f)=f$; 

$\bullet$  linearity on the linear span of $\ran(U_\mg) \cup \lin\{\delta_i\}_{i=1}^\infty$.

\noindent Then $W$ is bounded on $\ran(U_\mg)$, $L$ is not bounded on its domain, and $(L\delta_i)_{i=1}^\infty=(e_1, e_1, -e_1, e_2, e_1, -e_1, e_3, e_1, -e_1,\ldots)$, which is an s-pseudo-dual frame of $\mg$ and not a frame for $\h$.
\end{ex}

 \subsection*{Class of frames for which any s-pseudo-dual (resp. a-pseudo-dual) is necessarily a frame} \label{sdualdual}

As already noticed, for some frames $(g_i)_{i=1}^\infty$ there exist non-frame sequences $(f_i)_{i=1}^\infty$ satisfying 
  (\ref{fexp2}) and (\ref{fexp1}) (see e.g.  Example \ref{adn} with $p=2$). 
However, there exist frames  $(g_i)_{i=1}^\infty$ so that any sequence  $(f_i)_{i=1}^\infty$ satisfying   (\ref{fexp2}) and (\ref{fexp1}) is necessarily a frame. 
 Trivially, this is the case with every Riesz basis (as the only dual sequence of a Riesz basis is also a Riesz basis), but this may also happen with overcomplete frames. 
 Thus, it has been of interest in frame theory to characterize frames whose dual sequences are necessarily 
  frames. 
 The next statement gives a class of such frames.

\begin{prop} \label{qh1}
Let $\mg=\seqgr[g]$ be a frame for $\h$ and let 
 $\ker(T_\mg)$ be finite-dimensional\,\footnote{The condition {\it $\ker(T_\mg)$ being finite-dimensional} is equivalent to {\it $\mg$ being a near-Riesz basis} (i.e., to the existence of a finite set of elements of $\mg$, whose removal leaves a Riesz basis), see \cite{Holub}.}.  For a sequence $\seqgr[f]$ with elements from $\h$, the following statements are equivalent.
 \begin{itemize}
 \item[{\rm (i)}] $\seqgr[f]$ is an s-pseudo-dual of $\seqgr[g]$.
 \item[{\rm (ii)}] $\seqgr[f]$ is an a-pseudo-dual of $\seqgr[g]$.
 \item[{\rm (iii)}] $\seqgr[f]$ is a dual frame of $\seqgr[g]$.
  \end{itemize}
\end{prop}
\begin{proof}
That (iii) implies (i) and (ii) is obvious.

Assume that (ii) holds, i.e., that (\ref{fexp1}) holds. Since $\ker(T_\mg)$ is finite-dimensional, 
it follows from \cite{Holub} 
that the series $\sum_{i=1}^\infty c_i g_i$ converges only for $(c_i)_{i=1}^\infty\in\ell^2$. Therefore, 
$(\langle g,f_i\rangle)\in\ell^2$ for every $g\in\h$, which  implies that $\seqgr[f]$ is a Bessel sequence (see, e.g., \cite[Sec. 7.1]{Hbook}). Now \cite[Lemma 6.3.2]{Cbook} completes the proof that $\seqgr[f]$ is a dual frame of $\seqgr[g]$. 

Now assume that (i) holds. Let $g\in\h$. For every $f\in\h$, we have 
$\sum_{i=1}^\infty \langle f,g_i\rangle \langle f_i, g\rangle = \langle f,g\rangle$, 
so $\sum_{i=1}^\infty c_i \langle f_i, g\rangle $ converges for all $(c_i)_{i=1}^\infty\in \range(U_\mg)$. 
Furthermore, $\sum_{i=1}^\infty c_i \langle f_i, g\rangle $ converges  for every $(c_i)_{i=1}^\infty\in (\range(U_\mg))^\perp=\ker(T_\mg)$.  
Thus, $\sum_{i=1}^\infty c_i \langle f_i, g\rangle $ converges for all $(c_i)_{i=1}^\infty\in \ell^2$, which by \cite[Ex. 34.2]{Heuser} implies that $(\langle f_i, g\rangle)\in\ell^2$. 
Therefore, as above,  $\seqgr[f]$ is a Bessel sequence in $\h$ and hence a dual frame of $\seqgr[g]$. 
\end{proof}

\sloppy
As a simple illustration of Proposition \ref{qh1}, consider the frame $\seqgr[g] = (e_1, e_1, $ $e_2, e_3, e_4, e_5, \ldots)$ for $\h$. 
Any s-pseudo-dual (resp. a-pseudo-dual) of $\seqgr[g]$ has the structure $(w, e_1-w,$ $ e_2, e_3, e_4, e_5,  \ldots)$ for some $w\in\h$ and it is a frame for $\h$.

\subsection*{Acknowledgment}
The author is thankful to O. Christensen and H. Feichtinger for their valuable comments and suggestions.  
She also acknowledges support from the WWTF project MULAC  (\lq Frame Multipliers: Theory and Application in Acoustics\rq; MA07-025)  and the Austrian Science Fund (FWF)
START-project FLAME  (\lq\lq Frames and Linear Operators for Acoustical Modeling and Parameter Estimation\rq\rq; Y 551-N13).


\begin{thebibliography}{50}

\bibitem{Aldr} A. Aldroubi, 
\textit{Portraits of frames.}  Proc. Am. Math. Soc. {\bf 123} no.6 (1995), 1661--1668.

\bibitem{AST}
A. Aldroubi, Q. Sun,  W. Tang, \textit{$p$-frames and shift invariant subspaces of $L^p$.}
 J. Fourier Anal. Appl. {\bf  7} no.1 (2001),  1--21.

\bibitem{BLED} P. Balazs, B. Laback, G. Eckel, W.A. Deutsch, \textit{Time-frequency sparsity by removing perceptually irrelevant components using a simple model of simultaneous masking.} IEEE Trans. Audio Speech Lang. Processing {\bf 18} no.1 (2010), 34--49.

\bibitem{BSA}
 P. Balazs, D. T. Stoeva, J.-P. Antoine, \textit{Classification of general sequences by frame-related operators.}   Sampl. Theory Signal Image Process. {\bf  10} no.1--2 (2011), 151--170. 
 
\bibitem{Bary46}  
N. Bary, \textit{Sur les bases dans l’espace de Hilbert.}  
C. R. (Dokl.) Acad. Sci. URSS {\bf  54} (1946), 379--382. 
 
\bibitem{Bari}
  N. K. Bari, \textit{Biorthogonal systems and bases in Hilbert space}, Mathematics. Vol. IV, Uch. Zap. Mosk. Gos. Univ., {\bf 148}, Moscow Univ. Press, Moscow, 1951, 69--107. 

\bibitem{BL} M. Bownik, J. Lemvig, \textit{The canonical and alternate duals of a wavelet frame} Appl. Comput. Harmon. Anal. {\bf  23}  (2007), 263--272. 

\bibitem{Brezis}
 H. Brezis, \textit{Functional analysis, Sobolev spaces and partial differential equations}. 
Universitext. Springer, New York, 2011.


\bibitem{CasArt}
P. G. Casazza,  \textit{The art of frame theory.}  Taiwanese J. Math. {\bf  4} no.2 (2000),  129--201.



\bibitem{CCLL} P. Casazza, O. Christensen, S. Li, A. Lindner, 
  \textit{Riesz-Fischer sequences and lower frame bounds.}
  Z. Anal. Anwend. {\bf 21} no.2 (2002), 305--314.


\bibitem{CHL}
P. G. Casazza, D. Han,  D. R. Larson, \textit{Frames for Banach spaces.}
  Contemp. Math.  {\bf 247} (1999), The Functional and Harmonic Analysis of Wavelets and Frames, Baggett and Larson eds., 149--182.
 


\bibitem{Cinv} 
O. Christensen:  \textit{Frames and pseudo-inverses.} J. Math. Anal. Appl. {\bf  195} no.2 (1995), 401--414.

\bibitem{Cbook}
O. Christensen,   \textit{An Introduction to Frames and Riesz Bases.} Second Expanded Edition, 
 Series: Applied and Numerical Harmonic Analysis, 
Birkh\"auser, Boston, 
2016.




\bibitem{CCS}  P. G. Casazza, O. Christensen, D. T. Stoeva,
 \textit{Frame expansions in separable Banach spaces.} J. Math. Anal. Appl. {\bf 307} no.2 (2005), 710--723.

\bibitem{CS}
O. Christensen, D. T. Stoeva, \textit{$p$-frames in separable
Banach spaces.}   Adv. Comput. Math. {\bf 18} no.2-4 (2003), 117--126.


\bibitem{D90} I. Daubechies, \textit{The wavelet transform, time-frequency localization and signal analysis.}
 IEEE Trans. Inf. Theory {\bf 36} no.5 (1990), 961--1005.

\bibitem{Dbook} I. Daubechies,  \textit{Ten Lectures On Wavelets.} SIAM Philadelphia, 1992.



\bibitem{DGM} I. Daubechies, A. Grossmann, Y. Meyer, \textit{Painless nonorthogonal expansions.}
  J. Math. Phys. {\bf 27} (1986), 1271--1283.


\bibitem{DH} 
I. Daubechies, B. Han, \textit{The canonical dual frame of a wavelet
frame.}  Appl. Comput. Harmon. Anal. {\bf 12}   no.3 (2002),  269--285.



\bibitem{DMT}  P. Depalle, R. Kronland-Martinet, B. Torr\'{e}sani, \textit{Time-frequency multipliers for sound synthesis.} In: Proceedings of the Wavelet XII Conference, SPIE Annual Symposium, San Diego, 2007.

\bibitem{DSframe} R. J. Duffin,  A. C. Schaeffer, \textit{A class of nonharmonic Fourier
series.}   Trans. Am. Math. Soc. {\bf  72} (1952), 341-366.


\bibitem{FG}
 H. G. Feichtinger, K. Gr\"{o}chenig, \textit{Banach Spaces Related to Integrable Group
Representations and Their Atomic Decompositions I.}
  J. Funct. Anal. {\bf 86} no.2 (1989), 307--340.



\bibitem{FG2}
 H. G. Feichtinger, K. Gr\"{o}chenig,   \textit{Banach Spaces Related to Integrable Group
Representations and Their Atomic Decompositions II.}
  Monatsh. Math. {\bf 108} no.2-3 (1989), 129--148.


\bibitem{FW}
H. G. Feichtinger, T.Werther, \textit{Atomic systems for subspaces.}  Proceedings SampTA
2001 (L.Zayed, ed.), Orlando, 2001,  163--165.


\bibitem{FZim}
H. G. Feichtinger, G. Zimmermann,  \textit{A Banach space of test functions for Gabor analysis.} In: Gabor analysis and algorithms: Theory and applications, eds. H. G. Feichtinger and T. Strohmer, Birkh$\ddot{{\rm a}}$user, Boston 1998, 123--170. 


\bibitem{G}
K. Gr\"{o}chenig, \textit{Describing functions: atomic decompositions versus frames.}
 Monatsh. Math. {\bf  112} no.1 (1991), 1--42.

\bibitem{Gbook} 
K. Gr\"{o}chenig,  \textit{Foundations of Time-Frequency Analysis.} Series: Applied and Numerical Harmonic Analysis, Birkh\"auser, Boston, 2001.


\bibitem{HL}
D. Han, D. R. Larson, 
\textit{Frames, bases and group representations.}
Mem. Am. Math. Soc. {\bf  697}, American Mathematical Society (AMS), Providence, RI, 2000.



\bibitem{Hbook}
C. Heil, \textit{A Basis Theory Primer.} Expanded ed., Series: Applied and Numerical Harmonic Analysis, Birkh\"auser,  Basel, 2011. 
  
  \bibitem{Heuser}
 H. G. Heuser, \textit{Functional Analysis.  Transl. by John Horvath.} A Wiley-Interscience Publication. Chichester etc.: John Wiley \& Sons,  1982.

 \bibitem{Holub}
 J. R. Holub, \textit{Pre-frame operators, Besselian frames, and near-Riesz bases in Hilbert spaces.} 
 Proc. Amer. Math. Soc. {\bf 122} (1994), 779--785.

\bibitem{KA59}  L. V. Kantorovich, G. P. Akilov,  \textit{Functional Analysis in Normed Spaces.}
Gosudarstvennoe izdatelstvo Fiziko-Matematicheskoi Literatury, 
 Moskva, 1959 (in Russian). 

\bibitem{Li} S. Li,  \textit{On general frame decompositions.}
 Numer. Funct. Anal. Optimization,  {\bf 16} no.9-10 (1995),
1181--1191.

\bibitem{LO} S. Li, H. Ogawa, \textit{Pseudo-duals of frames with applications.} Appl. Comput. Harm. Anal. {\bf 11} no.2 (2001), 289-304.


\bibitem{MBKD} P. Majdak, P. Balazs, W. Kreuzer, M. D\"{o}rfler, \textit{A time-frequency method for increasing the signal-to-noise ratio in system identification with exponential sweeps.} In: Proceedings of the 36th IEEE International Conference on Acoustics, Speech, and Signal Processing, ICASSP 2011, 2011, pp. 3812--3815.

\bibitem{PS2}
S. Pilipovi\'c, D. T. Stoeva,  \textit{Series expansions in Fr\'echet spaces and their duals, construction of Fr\'echet frames.} J. Approx. Theory {\bf 163} (2011), 1729--1747. 

\bibitem{PS4}
S. Pilipovi\'c, D. T. Stoeva,  \textit{Fr\'echet frames, general definition and expansions.} Anal. Appl. {\bf 12} no.2 (2014), 195--208.


\bibitem{Rbook} W. Rudin,  \textit{Functional Analysis.} McGraw-Hill Series in Higher Mathematics, 
McGraw-Hill Book Comp., New York, 1973.

\bibitem{Slower} D. T. Stoeva,  \textit{Connection between the lower $p$-frame condition and existence of reconstruction formulas in a Banach space and its dual.} Ann. Sofia Univ., Fac. Math. Inf. {\bf 97} (2005), 123--133.

\bibitem{Srbasis} D. T. Stoeva, \textit{$X_d$-Riesz bases in separable Banach spaces.} \lq\lq Collection of papers, ded. to the 60th Anniv. of M. Konstantinov\rq\rq, 2008. 

\bibitem{Sb} 
D. T. Stoeva, \textit{Generalization of the frame operator and the canonical dual frame to Banach spaces.}  Asian-Eur. J. Math. {\bf 1} no.4 (2008), 631--643.

\bibitem{Sduals} D. T. Stoeva, \textit{$X_d$-frames in Banach spaces and their duals.}  Int. J. Pure Appl. Math. {\bf 52} no.1 (2009), 1--14.

\bibitem{Spert} D. T. Stoeva, \textit{Perturbation of frames for Banach spaces.}  Asian-Eur. J. Math.
{\bf 5} no.1 (2012), 1250011 (15 pages).

\bibitem{uncconv2011}
D. T. Stoeva, P. Balazs, \textit{Canonical forms of unconditionally convergent multipliers.} J. Math. Anal. Appl. {\bf 399} no.1 (2013), 252--259.

\bibitem{SBsampta15}
D. T. Stoeva, P. Balazs, \textit{The dual frame induced by an invertible frame multiplier.} 
In: Sampling Theory and Applications (SampTA), 2015 International Conference on, IEEE, 2015, 101--104. 


\bibitem{Young} R. M. Young, \textit{An {I}ntroduction to {N}onharmonic {F}ourier {S}eries.}  Academic Press,  New York, 1980. 

\end{thebibliography}
\end{document}